\documentclass{article}

\usepackage[english]{babel}
\usepackage[utf8x]{inputenc}
\usepackage[T1]{fontenc}

\usepackage{multirow}
\usepackage{mathrsfs}
\usepackage{amsfonts}
\usepackage{amsmath}
\usepackage{dsfont}
\usepackage{graphicx}
\usepackage{subfigure}%jia
\usepackage[colorinlistoftodos]{todonotes}
\usepackage[colorlinks=true, allcolors=blue]{hyperref}
\usepackage{amsthm}
\usepackage{float}
\usepackage{hyperref}
\usepackage{float}
\newtheorem{theorem}{Theorem}[section]

\newtheorem{definition}{Definition}[section]
\newtheorem{remark}[theorem]{Remark}

\numberwithin{equation}{section}
\numberwithin{lemma}{section}

\begin{document}
\title{Local estimators and Bayesian inverse problems with non-unique solutions}
\author{Jiguang Sun
\thanks{Department of Mathematical Sciences, Michigan Technological University, Houghton, MI 49931, U.S.A. ({\tt jiguangs@mtu.edu}).}
}
\date{}
\maketitle
\begin{abstract}
The Bayesian approach is effective for inverse problems. The posterior density distribution provides useful information of the unknowns. For problems with non-unique solutions, the classical estimators such as the maximum a posterior (MAP) and conditional mean (CM) are not suitable. We introduce two new estimators, the local maximum a posterior (LMAP) and local conditional mean (LCM). Some simple algorithms based on clustering to compute LMAP and LCM are proposed. Their applications are demonstrated by three inverse problems: an inverse spectral problem, an inverse source problem, and an inverse medium problem. 

\end{abstract}


\section{Introduction}
The Bayesian approach is an effective technique for inverse problems \cite{Fitzpatrick1991IP, KaipioSomerdalo2005, Stuart2010AN}. 
The problem is written in the form of statistical inferences.
Variables are viewed as being random and the known information can be coded in the priors. 
Using the Bayes' formula, one explores the posterior probability distribution of the unknowns.
This makes the Bayesian approach attractive for inverse problems with non-unique solutions, which often happens for partial data.

In practice, other than the posterior probability density, it is natural to provide some estimates of the unknowns, e.g., 
 maximum a posteriori (MAP) and conditional mean (CM), as the solutions for Bayesian inverse problems. However, such estimators might not carry sufficient information of the unknowns
 for complicate probability density, which motivates us to introduce new estimators. Having the inverse problems with non-unique solutions in mind,
 we introduce two new estimators, the local maximum a posterior (LMAP) and local conditional mean (LCM). They provide useful information for
 some probability densities and can be used to characterize the solutions of Bayesian inverse problems with non-uniqueness.
 
 It is relatively easy to compute the LMAPs for a given posterior density. Using the LMAPs, we propose an algorithm to compute LCMs based on
 the $k$-medoids or $k$-means for clustering problems. It works well for simple densities. Characterization of more complicate densities is challenging
 and needs further investigations.

To illustrate the applications of the LMAP and LCM, we consider three inverse problems: an inverse spectral problem, an inverse source problem, and an inverse medium problem.
All the problems have non-unique solutions given certain measurement data. Using the Bayesian inversion with the MCMC, we compute the posterior probability distributions
of the unknowns, which show clearly the existence of multiple solutions to the inverse problem. We employ the LMAP and LCM to characterize
the solutions (posterior probability densities) of these problems.

The rest of the paper is organized as follows.
In Section~\ref{BI}, we first give a brief introduction of the Bayesian inversion and the MCMC to explore the posterior density distribution. 
Then we define the LMAP and LCM. Section~\ref{Kmeans} describes how to compute the LMAPs and LCMs.
In Section~\ref{ISpec}, we consider an inverse spectral problem
to reconstruct the index of refraction given a Stekloff eigenvalue. 
% The posterior probability distribution indicates that there exist at least two solutions.
In Section~\ref{IMed}, the inverse problem is to reconstruct the wave speed using the data on a line segment.
In Section~\ref{ISou}, we consider a problem to reconstruct the location of an acoustic point source using the data at a single point.

\section{Bayesian Inversion and Local Estimators}\label{BI}
A simple statistical modal for the forward problem can be written as
\begin{equation*}
y=\mathcal{F}(x)+ \eta, \quad x\in X, \, y \in Y,
\end{equation*}
where $\mathcal{F}$ is the forward operator, $X, Y$ for some Banach spaces, and $\eta$ is the noise. 
For simplicity, we assume that the observation noise is normal with mean zero and independent of $x$, i.e., $\eta \sim {\mathcal N}(0,\sigma)$
with $\sigma$ being the covariance.

The inverse problem is to reconstruct $x$ given the noisy measurement $y$.
Denote by $\mu_0$ and $\mu_y$ the prior probability measure 
and the posterior probability measure of $x$, respectively. Let $\pi_0$ and $\pi_y$ denote the probability density functions of $\mu_0$ and $\mu_y$.
By Bayes' formula \cite{KaipioSomerdalo2005},
\begin{equation}\label{11}
\pi_y(x)=\frac{\pi_\eta\big(y-\mathcal{F}(x)\big)\pi_0(x)}
{\displaystyle\int_{X} \pi_\eta(y-\mathcal{F}(x))\pi_0(x)dx }.\\
\end{equation}
Thus
\[
\pi_y(x)\propto\pi_\eta\big(y-\mathcal{F}(x)\big)\pi_0(x)
\]
where $\propto$ means {\it proportional to}.  
The main task of Bayesian inversion is to explore the posterior density $\pi_y(x)$.

If the posterior density is obtained, point estimates such as maximum a posterior (MAP) or conditional mean (CM)
are often viewed as the solutions for the inverse problems. However, for complicate probability densities, 
the MAP and/or CM only provide partial information and, sometimes, can be misleading.
It is necessary to find ways to characterize such posterior probability densities. 

Motivated by the inverse problems with non-unique solutions,
we introduce two estimators: the local maximum a posterior (LMAP) and local conditional mean (LCM).
\begin{definition}\label{xLMAP} Denote the maximum a posteriori estimate of $\pi_y(x)$ by $x_{MAP}$, i.e.,
$x_{MAP} = \text{arg}\max_x \pi_y(x)$.
We call $x$ a local MAP, denoted by $x_{LMAP}$, if
\[
\pi_y(x) \ge \epsilon \, \max_x \pi_y(x) \quad \text{ and } \quad x = \text{arg}\max_{x \in N(x)} \pi_y(x)
\]
for some constant $\epsilon \in (0, 1)$ and $N(x)$, a neighborhood of $x$.
\end{definition}

\begin{definition}\label{xLCM} Denote by the conditional mean of $\pi_y(x)$ by $x_{CM}$, i.e., $x_{CM} = E\{x|y\} = \int_X x \pi_y(x) dx$. 
The local conditional mean $x_{LCM}$ is define as 
\[
x_{LCM} = \int_{S} x \pi_y(x) dx,
\]
where $S$ is a subset of $X$.
\end{definition}

Before we discuss how to compute $x_{LMAP}$ and $x_{LCM}$ and demonstrate their applications using three inverse problems with non-unique solutions,
we recall the MCMC (Markov chain Monte Carlo) to compute $\pi_{y}(x)$ in \eqref{11} \cite{KaipioSomerdalo2005}.

\vskip 0.2cm
{\bf MCMC Algorithm}\label{MHmethod}
\begin{itemize}
\item[1.] Pick the initial value $x_1$ and set $m\leftarrow 1$.
\item[2.] Draw $\tilde{x}$ from $\pi_0$ and calculate the acceptance ratio
	\[
		\alpha(x_k, \tilde{x}) = \min  \left\{1,  \frac{\pi(\tilde{x})}{\pi(x_k)}\right\}.
	\]
\item[3.] Draw $t \in [0, 1]$ from the uniform probability density. 
\item[4.] If $\alpha(x_n, \tilde{x}) \ge t$, set $x_{m+1} = \tilde{x}$, else $x_{m+1} = x_m$. 
	\begin{itemize}
		\item When $m=M$, the maximum number of iteration, stop.
		\item Otherwise, increase $m \leftarrow m+1$ and go to Step 2.
	\end{itemize}
\end{itemize}

\section{Characterization of $x_{LMAP}$ and $x_{LCM}$}\label{Kmeans}
% We propose an algorithm to compute $x_{LMAP}$ and $x_{LCM}$.
Assume a sampling of $\pi_y(x)$, still denoted by $\pi_y(x)$, is obtained. It is relatively straight forward to
decide $x_{LMAP}$. We first decide the value of $\epsilon$, e.g., $\epsilon = 0.5$, based on the knowledge of the inverse problem.
Then one may apply Def. \eqref{xLMAP} to find multiple $x^i_{LMAP}, i=1, \ldots, k$.

Given $k$ well-separated local MAPs $x^i_{LMAP}, i=1, \ldots, k$, it is reasonable to assume that the samples have $k$ clusters.
One can use $k$-medoids or $k$-means algorithm (see, e.g., \cite{CalvettiSomersalo2021}) to decide the clustering of the samples. 
We propose the following algorithm to find $x_{LMAP}$ and $x_{LCM}$ for $\pi_y(x)$.

\vskip 0.2cm
{\bf Algorithm LMAP-LMC}
\begin{itemize}
\item[] Given $\pi_y(x)$ and choose $\epsilon$.
\item[1.] Find $x_{LMAP}$s such that $x_{LMAP} > \epsilon x_{MAP}$. 
\item[2.] Assume there are $k$ well-separated $x^i_{LMAP}, i=1, \ldots, k$. Apply the $k$-medoids algorithm to find $k$ clusters of $\pi_y(x)$.
\item[3.] Exclude outliers in the clusters and compute $x^i_{LCM}, i=i, \ldots, k$.
\end{itemize}

\begin{remark}
The algorithm is effective for simple distributions (see Examples \ref{ISpec} and \ref{IMed}). For complicate posterior density distributions, 
additional knowledge of the inverse problems
and/or more powerful algorithms are needed (see Example \ref{ISou}). How to characterize the posterior density function $\pi_y(x)$ is 
an interesting and important topic.
\end{remark}

\begin{remark}
The number of clusterings $k$ can be obtained without using $x^i_{LMAP}$'s.
For example, one can use the within-cluster mean distance to decide $k$ (see, e.g., Section 3.4 of \cite{CalvettiSomersalo2021}). 
\end{remark}

\begin{remark}
The exclusion of the outliers in the clusters is consistent with the definitions of the local estimators. 
\end{remark}

\section{Inverse Problems with Non-unique Solutions}
We consider three inverse problems with non-unique solutions: an inverse spectral problem,
an inverse medium problem, and an inverse source problem. In fact, the last example shows that, for complicate posterior density functions,
even the local point estimators might not be enough. More advanced estimators such as curve estimators or set estimators should be considered.

\subsection{Inverse Spectral Problem}\label{ISpec}
Let $D \subset \mathbb R^2$ be a disk with radius $1$ with boundary $\partial D$. Let $k$ be the wavenumber and $n$ be the index of refraction, which is a real constant.
We call $\lambda$ a Stekloff eigenvalue if there exists a non-trivial function $w$ such that
\begin{equation}\label{stekloff}
  \left\{
   \begin{array}{rl}
   \Delta w+\beta^2n w=0, &\textrm{in}\  D, \\
   \partial w/\partial \nu+\lambda w=0, &\textrm{on}\ \partial D.
   \end{array}\right.
\end{equation}
When $\beta^2$ is not a Dirichlet eigenvalue of $D$ (see Chp. 3 of \cite{SunZhou2016}), 
the Stekloff eigenvalues are real and discrete.
For simplicity, we assume that $\beta=1$ and $n$ is a real constant. Consider the following inverse problem:
\begin{itemize}
\item[{\bf IP1}] Given a Stekloff eigenvalue $\lambda$, find the index of refraction $n$.
\end{itemize}
The above problem can be written as
a statistical inference for $n$ \cite{LiuLiuSun2019IP}
~\begin{equation} \label{eq 3.21}
\lambda=\mathcal{F}(n)+\eta,
\end{equation}
where $\lambda$ is the given Stekloff eigenvalue, $n$ is the unknown random variable, $\mathcal{F}$ is the
operator mapping $n$ to $\lambda$ based on the partial differential equation \eqref{stekloff}. 
We assume that an a priori that $a < n < b$ where $a, b$ are two real constants.
A natural choice for the prior is $n \sim \mathcal{U}(a,b)$, 
where $\mathcal{U}$ denotes the uniform distribution. 

% The goal of the Bayesian inverse problem is to explore the conditional probability distribution
% $\pi_{\lambda}(n)$, called the posterior distribution of $n$.
By the Bayes' formula, the posterior distribution satisfies
\begin{equation} \label{eq 3.22}
\pi_{\lambda}(n)\propto \mathcal{N}({\lambda}-\mathcal{F}(n),\,\sigma^{2}) \times \mathcal{U}(a,b),
\end{equation}
i.e.,
\begin{equation} \label{pipostnlambda}
\pi_{\lambda}(n) \propto \exp \Big(-\frac{1}{2 \sigma^{2} }|{\lambda}-\mathcal{F}(n)|\Big)\times I(a\leq n\leq b),
\end{equation}
where $I$ is the density function for $\mathcal{U}(a,b)$.
% The inverse problem is to   the posterior distribution $\pi(n|{\boldsymbol \lambda})$.

\begin{figure}[h!]
\begin{center}
{ \scalebox{0.50} {\includegraphics{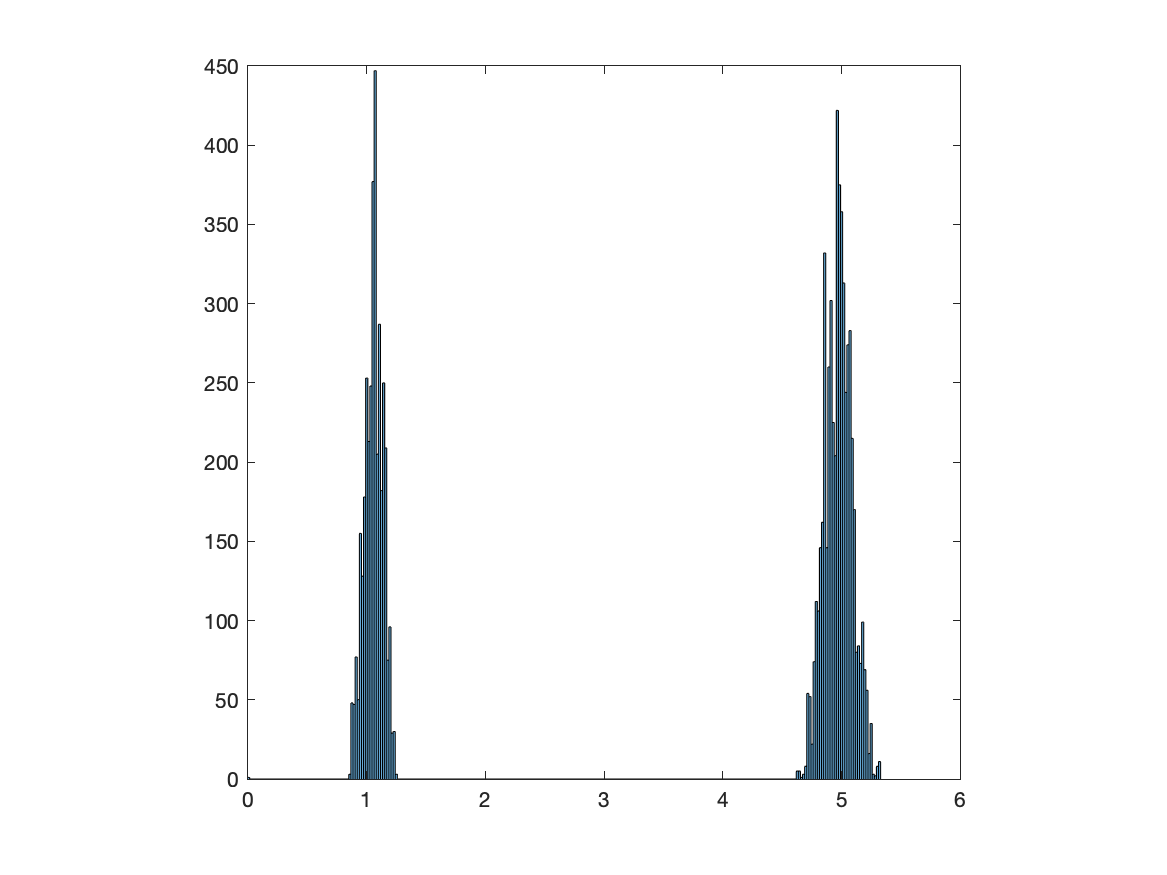}}}
\caption{Histogram of the posterior probability distribution for $n$.}
 \label{Fig1}
\end{center}
\end{figure}

Let $\lambda = 0.62$ be a given Stekloff eigenvalue \cite{LiuLiuSun2019IP}. Assume that $n \sim \mathcal{U}(0, 6)$.
Set $\sigma = 0.05$ and $K=10000$ and carry out the MCMC algorithm. At each iteration,
we compute the eigenvalue of \eqref{stekloff} closest to $\lambda$ using a finite element method \cite{SunZhou2016}.
We discard the first $1000$ samples and show the histogram in Fig.~\ref{Fig1}. 

For the posterior density,  $n_{CM} = 3.413$ and $n_{MAP} = 1.075$.
For $n=n_{MAP}$, the eigenvalue closest to $0.62$ is $\lambda^*=0.625$.
However,  it is clear from Fig.~\ref{Fig1} that $n_{CM}$ does not provide a reasonable answer to the inverse problem.
In fact, for $n=n_{CM}$, the Stekloff eigenvalue closest to  $0.62$ is $\lambda^* = 0.005$. 
The posterior density function has two local maximums (the samples have two clusters). 
In Table~\ref{Table1}, we set $\epsilon = 0.5$ and show the local estimators computed using {\bf LMAP-LCM}.
The associated eigenvalues $\lambda^*$ are listed as well, which are in good accordance with $\lambda=0.62$. 
\begin{table}[h!]
\centering
\begin{tabular}{r|r|r|r }
\hline
$n^1_{LMAP} = 1.075$& $n^2_{LMAP} = 5.005$ &$n^1_{LCM} = 1.063$&$n^2_{LCM} = 4.973$ \\\hline
$\lambda^*=0.625$ & $\lambda^* =0.616$ & $\lambda^*=0.648$ & $\lambda^*=0.614$ \\\hline
\end{tabular}
\label{Table1}
\caption{Local estimators for $n$ and the associated eigenvalue $\lambda^*$ closet to $\lambda$.}
\end{table}

\subsection{Inverse Medium Problem}\label{IMed}
We consider an inverse medium problem from \cite{Ramm1995}.
Let $D =(0, 1)^2$ and assume that $S \subset \partial D$ is on the plane $x_2=0$.
Let $f(x,t)$ be a source and $c(x)>0$ be the velocity.
The acoustic pressure $u$ solves the problem
\begin{subequations}\label{IMP}
\begin{align}
\label{IMPE} c^{-2}(x) u_{tt}-\triangle u &=f(x, t) \qquad &\text{in } D \times [0, \infty),\\[1mm]
\label{IMPB}\partial u/\partial \nu&=0 \qquad &\text{on } \partial D,\\[1mm]
\label{IMPI} u =  u_t&=0 \qquad &\text{at } t=0.
\end{align}
\end{subequations}

Assuming that $c(x)=c$, we consider the following inverse medium problem.

\begin{itemize}
\item[{\bf IP2}] Given the data $u(x,t), x \in S, t > 0$, recover the constant speed $c$.
\end{itemize}

The solution of \eqref{IMP} is given by
\begin{equation}\label{exactu}
u(x, t) = \sum_{m=0}^\infty =  u_m(t)\phi_m(t), \quad m=(m_1, m_2), 
\end{equation}
where $\phi_m = \gamma_{m_1m_2} \cos(m_1x_1)\cos(m_2 x_2)$ and
\begin{eqnarray*}
&&\int_D \phi_m^2(x) dx =1, \quad \triangle \phi_m + \lambda_m \phi_m = 0, \quad \lambda_m:=m_1^2+m_2^2, \\
% &&\partial \phi_n/\partial \nu = 0 \quad \text{ on  } \partial D, \quad \\
&& \gamma_{00} = 1/\pi, \quad \gamma_{m_1 0}  = \gamma_{0m_2} = \sqrt{2}/\pi, \quad \gamma_{m_1m_2} =  2/\pi \text{ if } m_1, m_2 > 0,\\
&& u_m(t):=u_m(t,c) = \frac{c}{\sqrt{\lambda_m}} \int_0^t \sin \left(c \sqrt{\lambda_m}(t-\tau) \right) f_m(\tau)d\tau, \\
&&  f_m(t) = \int_D  f(x,t)\phi_m(x) dx.
\end{eqnarray*}

The solution $u$ on $S$ are thus given by
\begin{equation}\label{udata}
u(x_1, 0, t)  = \sum_{m=0}^\infty u_m(t,c) \gamma_{m_1 m_2} \cos(m_1 x_1).
\end{equation}
%The data $u(x_1, 0, t)$ for $c=c_1$ and $c=c_2$ are the same if and only if
%\begin{equation}\label{gammam1m2}
%\sum_{m_2}^\infty \gamma_{m_1 m_2} u_m(t,  c_1) = \sum_{m_2}^\infty \gamma_{m_1 m_2} u_m(t,  c_1) , \quad \forall t > 0, \quad  \forall m_1 \ge 0.
%\end{equation}
%Taking the Laplace transform of \eqref{gammam1m2}, one obtains
%\begin{equation}\label{overlinefmp}
%\sum_{m_2=0}^\infty \gamma_{m_1 m_2} \overline{f}_m(p) \left( \frac{c_1^2}{p^2+c_1^2\lambda_m} - \frac{c_2^2}{p^2+c_2^2\lambda_m} \right) = 0, 
%\quad \forall p > 0, \forall m_1.
%\end{equation}

%For $c_1\ne c_2, c_1, c_2 > 0$, there exists infinitely many ways to find $\overline{f}_m(p)$ such that \eqref{overlinefmp} holds.
%For example, let $c_1=1$, $c_2  = 2$, $\overline{f}_{m_1m_2} = 0$ for $m_1 \ne 0$, $m_2 \ne 1$, or $m_2 \ne 2$,
%$\overline{f}_{01}(p) = -(p^2+1)/((p+1)(p^2+16))$, $\overline{f}_{02}(p) = 1/(p+1)$. Then \eqref{overlinefmp} holds. 
%Therefore, if
%\[
%f(x, t) = \frac{\sqrt{2}}{\pi} \left( f_{01}(t)\cos(x_2)+f_{02}(t)\cos(2x_2)\right), c_1=1, c_2=2,
%\]
%where
%\[
% f_{01}(t) = -\frac{2}{17} \exp(-t) - \frac{15}{17}\left( \cos(4t)-\frac{1}{4}\sin(4t) \right), \quad  f_{02}(t) = - \exp(-t).
%\]
Let $u_1(x, t, c)$ be the solution of \eqref{IMP} for some speed $c$. 
It was shown in \cite{Ramm1995} that $u_1(x, t, 1)=u_2(x, t, 2)$ for $x \in S$ and $t >0$.

%\begin{figure}[h!]
%\begin{center}
%{ \scalebox{0.50} {\includegraphics{IMP.eps}}}
%\caption{Histogram of the posterior probability distribution for $c$.}
% \label{Fig2}
%\end{center}
%\end{figure}
For the inverse problem, the data are given by ${\boldsymbol u}:=\{u^m_i\}_{i=1}^{10}$, $u_i = u((i+1)h-h/2, 0, 1)$, $h=1/10$. 
The measured data are obtained by \eqref{exactu} with $c=1$ adding 5\% of white noise. The statistical model for the problem is 
\begin{equation} \label{uFeta}
{\boldsymbol u}=\mathcal{F}(c)+\eta,
\end{equation}
where $\mathcal{F}$ is the forward operator by \eqref{IMP} and $\eta$ is the noise.

\begin{figure}[h!]
\begin{center}
\begin{tabular}{lll}
\resizebox{0.5\textwidth}{!}{\includegraphics{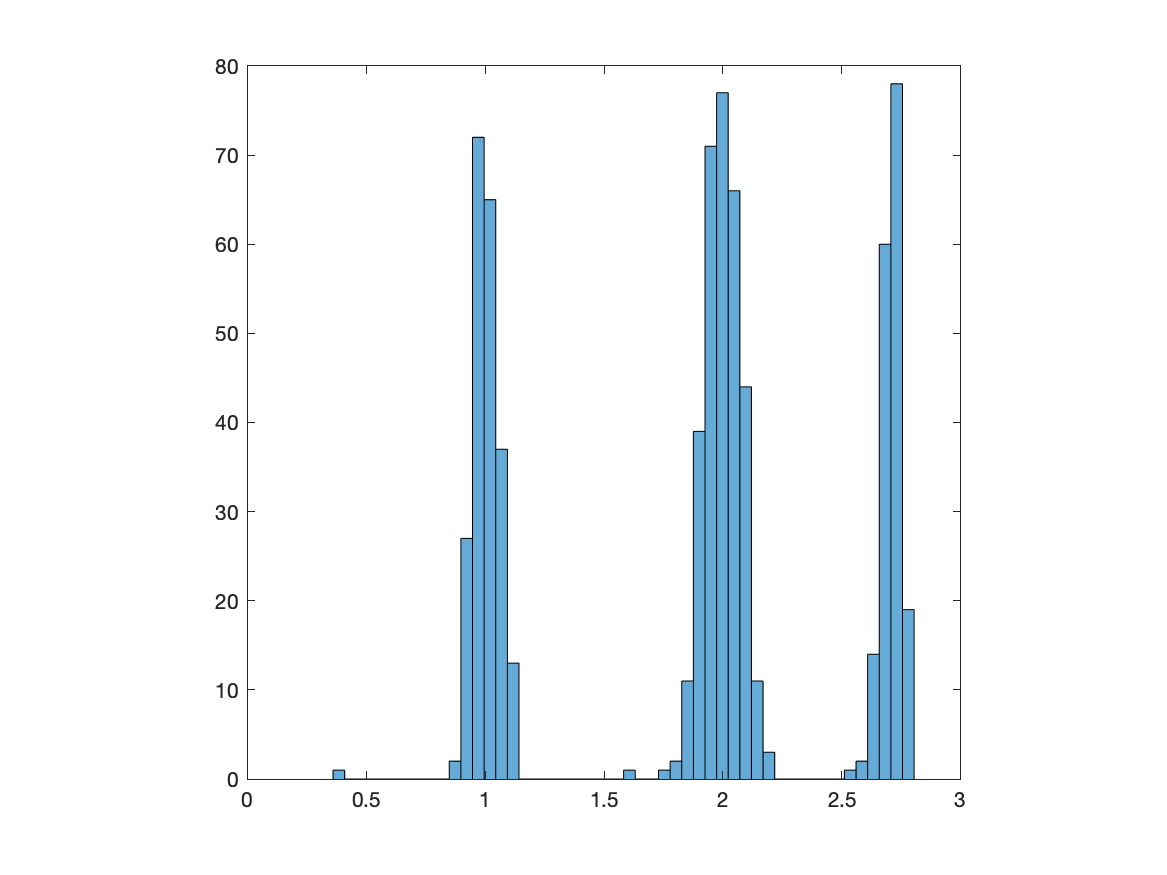}}&
\resizebox{0.5\textwidth}{!}{\includegraphics{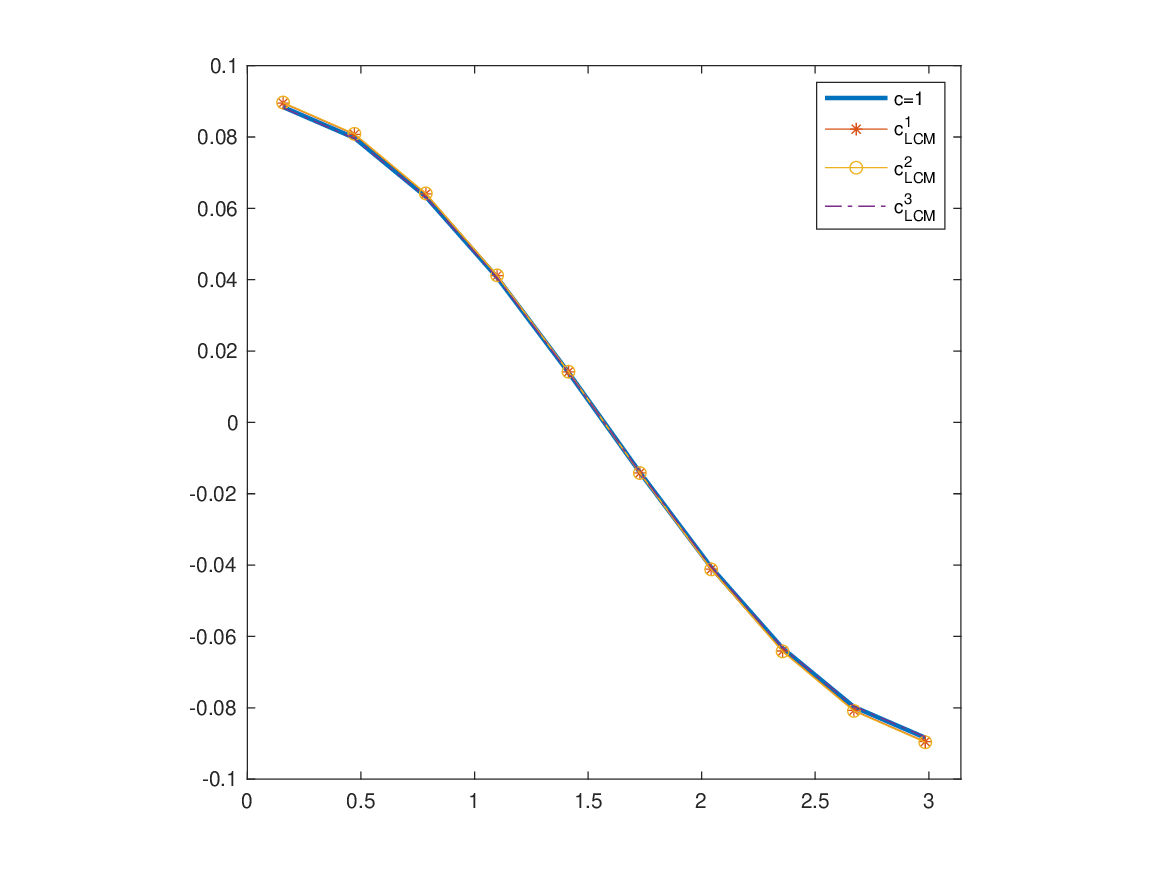}}
\end{tabular}
\end{center}
\caption{Left: histogram. Right: Values of $u$ using LCMs.}
\label{Fig2}
\end{figure}

Letting the prior density function for $c$ to be $\mathcal{U}(0, 3)$, we employ the MCMC algorithm to compute
the posterior density function using 5000 samples. The histogram is shown in the left picture in Fig.~\ref{Fig2}.
As we expected, there are many samples accumulate around $1$. In addition, there are also
large number of samples around $2$ and $2.7$. We compute the local estimators and show them in Table~\ref{Table2}.
In the right picture of Fig.~\ref{Fig2}, we show the exact values of $u (c=1)$ and the values by the three LCMs for $c$.
They coincide very well.

\begin{table}[h!]
\centering
\begin{tabular}{r|r|r }
\hline
$c^1_{LCM} = 1.011$& $c^2_{LCM} = 1.984$ &$c^3_{LCM} = 2.713$\\ \hline
$c^1_{LMAP} = 0.975$& $c^2_{LMAP} = 1.985$ &$c^3_{LMAP} = 2.705$\\ \hline
\end{tabular}
\label{Table2}
\caption{Local estimators for $c$.}
\end{table}

%%%%%%%%%%%%%%%%%%%%%%%%%%%%%%%%%%%%%%%%%%%%%%%%%%%%%%%%%%%%%%%%%%%%%%%
\subsection{Inverse Source Problem}\label{ISou}
We present an inverse source problem to show that it is challenging to characterize complicated posterior density functions.
%Let $F\in L^2(\mathbb{R}^2)$ denote the source with $\text{supp}\, F \subset V$, 
%where $V$ is a bounded domain on $\mathbb{R}^2$. 
Consider the time-harmonic acoustic wave field $u\in H^1_{loc}(\mathbb{R}^2)$ radiated by a point source at $z$ such that
\begin{subequations}\label{Helmholtz1}
\begin{equation}\label{Helmholtz2}
\Delta u+\beta^2 u= -  \delta(x-z) \quad \text{in}\; \mathbb{R}^2,
\end{equation}
\begin{equation}\label{Sommerfeld}
\lim\limits_{r \rightarrow \infty} \sqrt{r}
\left(\frac{\partial{u}}{\partial r}-i \beta u\right)=0, \quad r=|x|,
\end{equation}
\end{subequations}
where $\beta$ is called the wavenumber, \eqref{Helmholtz2} is the Helmholtz equation and \eqref{Sommerfeld} is the Sommerfeld radiation condition.
The solution to \eqref{Helmholtz1} is given by
\begin{equation}\label{uxzk}
u(x, z, \beta) :=
\dfrac{i}{4}H^{(1)}_0(\beta|x-z|),
\end{equation}
where $H_0^{(1)}$ is the Hankel function of zeroth order and first kind \cite{ColtonKress2013}.

The inverse source problem (ISP) is to determine the location $z$ of the point source from the measurement of $u$ at a point $x_0$.

\begin{itemize}
\item[{\bf IP3}] Given the data $u(x_0, z, 1), x_0=(0, 3)$, find the source location $z$.
\end{itemize}

We write the statistical model as
~\begin{equation} \label{SISP}
u=\mathcal{F}(z)+\eta,
\end{equation}
where $\mathcal{F}$ is the forward operator given by \eqref{uxzk} and $\eta$ is the noise given by the normal distribution $\mathcal{N}(0, 0.002)$.
The prior for $z$ is the uniform distribution $\mathcal{U}([-2,2]\times [-2, 2])$. 

\begin{figure}[h!]
\begin{center}
\begin{tabular}{lll}
\resizebox{0.5\textwidth}{!}{\includegraphics{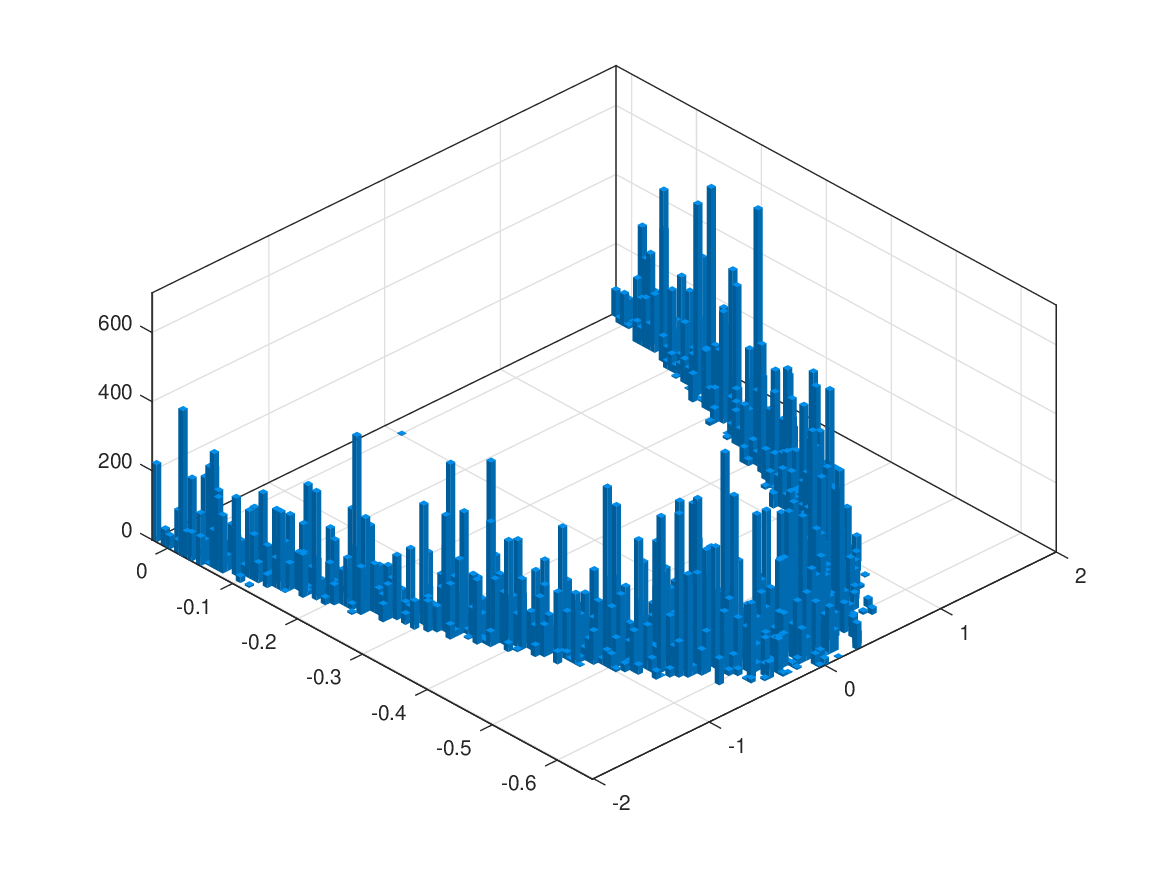}}&
\resizebox{0.5\textwidth}{!}{\includegraphics{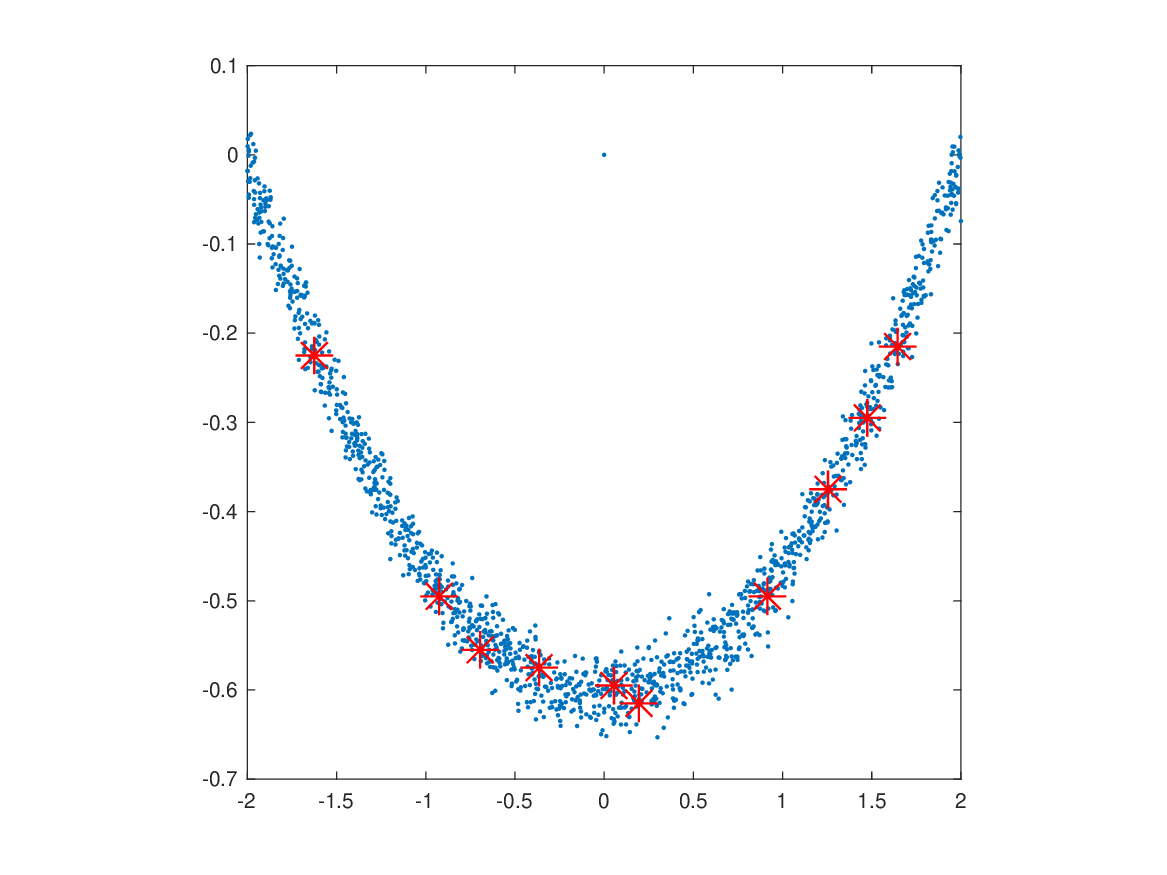}}
\end{tabular}
\end{center}
\caption{Left: histogram. Right: accepted samples. '*' indicate the LMAPs.}
\label{Fig3}
\end{figure}

The given value $u(x_0, z_e, \beta)$ is computed using \eqref{uxzk} with $z_e=(2, 0)$ and $\beta=1$, adding 5\% uniformly distributed noise.
In the MCMC, we draw 100,000 samples. The histogram and accepted samples are shown in Fig.~\ref{Fig3}. 

Any $z$ on the circle $C$ centered at $x_0=(0, 3)$ with radius $r=\sqrt{3^2+2^2}$ would give the same value as $u(x_0, z_e, 1)$.
Indeed, the samples accumulate around the curve $C \cap [-2, 2]^2$. 
The LMAPs provide the possible solutions but are only partial (see the right picture of Fig.~\ref{Fig3}).
This example shows that, for complicate posterior density functions, how to identify $S$ and define the local estimators
are problem-dependent and challenging. 

\end{document}